\documentclass[11pt,titlepage,a4paper]{article}
\usepackage{bozhomac}	
\usepackage{bozhlogo}	
\usepackage{cite}	
\usepackage{graphicx}	
\usepackage{floatflt}	

\title{\bfseries	\vspace*{-1.678902345in}
{\huge Linear vector fields and\\[1ex] exponential law}
}

\vspace{1.7ex}

\author{
Bozhidar Z. Iliev
\thanks{Laboratory of Mathematical Modeling in Physics,
Institute for Nuclear Research and \mbox{Nuclear} Energy,
Bulgarian Academy of Sciences,
Boul.\ Tzarigradsko chauss\'ee~72, 1784 Sofia, Bulgaria}
\thanks{E-mail address: bozho@inrne.bas.bg}
\thanks{URL: http://theo.inrne.bas.bg/$\sim$bozho/}
\\
Maido Rahula
\thanks{Institute of pure mathematics, Faculty of mathematics and informatics,
University of Tartu, J.~Liivi~2, Tartu 50409, Estonia}
\thanks{E-mail address: rahula@math.ut.ee}
\thanks{URL: http://math.ut.ee/pmi/inimesed/rahulaeng.html}
}

\date{
 \vspace{2.27ex}\ShortTitle{Linear vector fields and exponential law}\\[0.27ex]
 \vspace{3.27ex}
\small
	\begin{tabular}{r@{$\colon\to~$}l}
 \vspace{0.09ex} Last update	& April 1, 2006	\\[0.09ex]
 \vspace{0.27ex} Produced	& \fbox{\today}	\\[0.27ex]
	\end{tabular} \\[1.27ex]
 \small
	\begin{tabular}{r@{$\colon~$}l}
 \normalsize\sffamily\bfseries
  \vspace{0.27ex} http://www.arXiv.org e-Print archive No. &
 \normalsize\sffamily\bfseries
 math.DG/0604005	\\[1.27ex]
	\end{tabular} \\[-0.27ex]
\normalsize
 \vspace{4.27ex}{\Huge\BOZHO}	\\[4.27ex]
	\begin{tabular}{r@{\hspace{0.512em}}|@{\hspace{0.512em}}l}
 \vspace{0.27ex}\MSC[2001]{53B99, 38A17, 58J70\\ 34A26, 34A99}
&
 \vspace{0.27ex}\PACS[2001]{02.30.Hq, 02.30.Jr\\ 02.40.Mq, 02.40.Vh}
	\end{tabular} \\[1.27ex]
\vspace{0.27ex}\KeyWords{Vector and Tensor fields, Linear vector fields\\
	Fundamental vector fields, Exponential law\\
	Flows of vector fields, Classification of flows\\
	Representations of the general linear group}	\\[0.27ex]
}

\begin{document}		

\renewcommand{\thepage}{\roman{page}}

\renewcommand{\thefootnote}{\fnsymbol{footnote}} 
\maketitle				
\renewcommand{\thefootnote}{\arabic{footnote}}   

\tableofcontents		


\begin{abstract}

	The paper is devoted to vector fields on the spaces $\field[R]^2$ and
$\field[R]^3$, their flow and invariants. Attention is plaid on the tensor
representations of the group $\GL(2,\field[R])$ and on fundamental vector
fields. The rotation group on $\field[R]^3$ is generalized to rotation groups
with arbitrary quadrics as orbits.

\end{abstract}

\renewcommand{\thepage}{\arabic{page}}


\section {Introduction}
\label{Introduction}

	In the paper are presented equations for calculation of invariants of a
vector field when an infinitesimal symmetry is known. The flow and invariants
of a linear vector field, which is a fundamental vector field (group operator)
for the general linear group $\GL(2,\field[R])$ are governed by exponential
law. Representations of $\GL(2,\field[R])$ on tensor fields are investigated. If
a tensor field is invariant relative to a vector field, on its components act
an extension of this field (classical situation); if the components are
invariant, the field and its Lie derivatives form linear ordinary differential
equation. Some examples are presented.

	The flow and invariants of vector fields on $\field[R]^2$ are discussed
in section~\ref{Sect0}. Linear vector fields on $\field[R]^2$, their Lie
derivatives and flows are briefly considered in section~\ref{Sect0-1}.
	Section~\ref{Sect9} contains a full classification of the  flows of the
linear vector fields on $\field[R]^2$. In particular, it is recalled that their
singularities can of the types focus, saddle point and node; the corresponding
phase portraits are drawn. The basic fundamental vector fields of the rotation
group of $\field[R]^3$ are investigated in section~\ref{Sect10}.
Section~\ref{Sect11} is devoted to tensor representations of the two\ndash
dimensional real general linear group.

	Section~\ref{Conclusion} closes the paper.


\section{Flow and invariants of vector fields on $\field[R]^2$}
	\label{Sect0}%

	Let on $\field[R]^2$, coodinated by $(u,v)$, be given a vector field
    \begin{equation}    \label{0.1}
X=x\frac{\pd }{\pd u} + y\frac{\pd }{\pd v}
    \end{equation}
with components $(x,y)$ depending on $u$ and $v$. The action of $X$ on a $C^1$
function $f$ will be denoted via ``prime'', \viz
\[
f'=f_1x+f_2y
\]
with $f_1:=\frac{\pd f}{\pd u}$ and $f_2:=\frac{\pd f}{\pd v}$. With a prime
will be denoted also the Lie derivatives relative to $X$ of tensor fields, if
the field $X$ is clear from the context; otherwise the standard notation
$\Lied_X$ will be used.
	In particular,  we have
    \begin{equation}    \label{0.2}
    \begin{split}
u'&=x(u,v) \\
v'&=y(u,v) .
    \end{split}
    \end{equation}

	If we consider these equalities as equations, they define the flow of $X$, \viz
a 1\ndash parameter group $a_t$, $t\in\field[R]$, of diffeomorphisms of
$\field[R]$ such that
    \begin{equation}    \label{0.3}
a_t \colon (u,v)\to (u_t,v_t).
    \end{equation}
A function $f$ is dragged by the flow $a_t$ according to the rule
    \begin{equation}    \label{0.4}
f \mapsto  f_t=f\circ a_t.
    \end{equation}
It it happens that $f_t$ has a $C^\omega$ dependence on $t$, we can write the
Maclaurin series
\[
f_t = \sum_{k=0}^{\infty} f^{(k)} \frac{t^k}{k!} .
\]
	If between the derivatives of  $f$ there is some connection, \ie an
ordinary differential equation relative to $f$, its solution gives some dragged
function; \eg it $f''+f=0$, then $f_t=f\cos t + f'\sin t$.

	The differential form
    \begin{equation}    \label{0.5}
\omega=-y \od u + x \od v ,
    \end{equation}
which is annihilated by $X$, is a useful tool for finding invariants of $X$. If
$\omega$ is closed, \ie $\od \omega=0$, it is locally exact and equals the
differential of an invariant $I$ of $X$:
    \begin{equation}    \label{0.6}
\omega=\od I \ \Rightarrow\ I'=\od I(X)=\omega(X)=0,
    \end{equation}
so that
    \begin{equation}    \label{0.7}
I=\int\omega.
    \end{equation}
If $\omega$ is not exact, there is an integrating factor $\mu$ such that
$\mu\omega$ is an exact form, $0=\od(\mu\omega)=(\mu'+\diver X\cdot\mu)\omega$
with $\diver X=x_1+y_2$. Consequently, solution of the equation
    \begin{equation}    \label{0.8}
\mu'+\diver X\cdot\mu = 0
    \end{equation}
is an integrating factor of $\omega$ and hence the equation
    \begin{equation}    \label{0.9}
I=\int\mu\omega
    \end{equation}
gives an invariant of $X$.

	If a vector field $P$ is an infinitesimal symmetry of $X$, \ie
\[
\Lied_PX\parallel X \text{ or } \Lied_P\omega\parallel \omega ,
\]
where $\parallel$ means equal up to multiplicative function,
then $\frac{1}{\omega(P)}$ is an integrating factor for $\omega$
(see~\eref{0.8}) and hence
    \begin{equation}    \label{0.10}
I=\int\frac{\omega}{\omega(P)}
    \end{equation}
is an invariant of $X$. This fact is a consequence of $\omega(X)=0$ and
$\omega'=\diver X\cdot\omega$ which imply $(\omega(P))'=\diver
X\cdot\omega(P)$.

	In particular, if the coordinate functions $u$ and  $v$ are
respectively a canonical parameter and invariant of $P$,  $P(u)=1$ and
$P(v)=0$, then we can set $P=\frac{\pd }{\pd u}$, so that $\Lied_PX\parallel X$
and $\omega(P)=-y=-X(v)$. Therefore
 $\frac{\omega}{\omega(P)}=\od u-f(v)\od v$ with $f(v):=\frac{X(u)}{X(v)}$,
due to $x=X(u)$, and consequently $f(v)$ is an invariant of $P$
and
    \begin{equation}    \label{0.11}
I=u - \int f(v)\od v
    \end{equation}
is an invariant of $X$.


\section{Linear vector fields on \protect{$\field[R]^2$}}
\label{Sect0-1}

	A linear vector field on $\field[R]^2$ coordinated by $(u,v)$ is of the
form
    \begin{equation}    \label{0.12}
X=u'\frac{\pd }{\pd u} + v'\frac{\pd }{\pd v}
    \end{equation}
where its components $u'$ and $v'$ are linear homogeneous functions of the
coordinates $u$ and $v$, \ie
    \begin{equation}    \label{0.13}
\begin{pmatrix} u' \\ v' \end{pmatrix}
=
    \begin{pmatrix} c_1 & c_2 \\[0.3ex]
		    c_3 & c_4
    \end{pmatrix}
\begin{pmatrix} u \\ v \end{pmatrix} .
    \end{equation}
Introducing the matrices
    \begin{equation}    \label{3.2-1}
    \begin{split}
U := \begin{pmatrix} u \\ v \end{pmatrix} \quad
R := \frac{\pd }{\pd U}
  :=  \Bigl( \frac{\pd }{\pd u} \, , \frac{\pd }{\pd v} \Bigr)
\quad
C:=
    \begin{pmatrix} c_1 & c_2 \\[0.3ex]
		    c_3 & c_4
    \end{pmatrix}  ,         \quad
\theta:= \od U
       := \begin{pmatrix} \od u \\ \od v \end{pmatrix}
    \end{split}
    \end{equation}
we see that the flow of $X$ is locally governed by an exponential law:
    \begin{equation}    \label{0.14}
U'=CU \ \Rightarrow\ U_t=\e^{tC} U.
    \end{equation}
	One can interpret this as follows: from the l.h.s.\ of the implication
is an ordinary differential equation and from its r.h.s.\ is written its
general solution.

	Calculating the Lie derivatives of the frame $\frac{\pd }{\pd U}$ and
coframe $\od U$ and their dragging in the flow of X, we get
    \begin{subequations}    \label{0.15}
    \begin{align}    \label{0.15a}
\Bigl(\frac{\pd }{\pd U}\Bigr)' = - \frac{\pd }{\pd U} C
&\ \Rightarrow\
\Bigl(\frac{\pd }{\pd U}\Bigr)_t =  \frac{\pd }{\pd U} \e^{-tC}
\\		    \label{0.15b}
(\od U)' = C U
&\ \Rightarrow\
(\od U)_t = \e^{tC} \od U .
    \end{align}
    \end{subequations}

	From the Hamilton-Cayley formula for $C$, \viz
 $C^2-\tr C\cdot C+\det C\cdot\openone$=0, and the l.h.s.\ of~\eref{0.14}, we
see that (each element of) $U$ satisfies the second-order ordinary differential
equation
    \begin{equation}    \label{0.16-0}
U''-\tr C\cdot U' + \det C\cdot U = 0 .
    \end{equation}
Similar equations for $\od U$ and $\frac{\pd }{\pd U}$ can be found by means of
the right hand sides of~\eref{0.15}:
    \begin{equation}    \label{0.16-1}
\Bigl(\frac{\pd }{\pd U}\Bigr)'' +
	\tr C\cdot \Bigl(\frac{\pd }{\pd U}\Bigr)' +
	\det C\cdot \Bigl(\frac{\pd }{\pd U}\Bigr) = 0
    \end{equation}
    \begin{equation}    \label{0.16-2}
(\od U)''-\tr C\cdot (\od U)' + \det C\cdot \od U = 0 .
    \end{equation}

	If $Y$ and $\Phi$ are respectively a vector field and one-form with
matrix representations $Y=\frac{\pd }{\pd U}y$ and $\Phi=\varphi\od U$, with
$y=\begin{pmatrix} y^1 \\ y^1 \end{pmatrix}$ and
$\varphi=(\varphi_1,\varphi_2)$, one can easily deduce the implications
    \begin{subequations}    \label{0.16}
    \begin{align}    \label{0.16a}
Y'=0
&\ \Rightarrow\
y'=Cy
&\ \Rightarrow\
y''-\tr C\cdot y' + \det C\cdot y = 0
&\ \Rightarrow\
y_t = \e^{Ct} y
\\		    \label{0.16b}
\Phi'=0
&\ \Rightarrow\
\varphi'= - \varphi C
&\ \Rightarrow\
\varphi''  + \tr C\cdot \varphi' + \det C\cdot \varphi = 0
&\ \Rightarrow\
\varphi_t = \varphi \e^{-Ct} .
    \end{align}
    \end{subequations}

	We see here a representation of the group $\GL(2,\field[R])$ on the
space of components $y$ and/or $\varphi$. To an element $C\in\gl(2,\field[R])$
(here we identify $\gl(2,\field[R])$ with the isomorphic to it set of all
$2\times2$ matrices) corresponds the one-parameter group $\e^{Ct}$ in
$\GL(2,\field[R])$, which defines the dragging of $y$ and $\varphi$ in the flow
of $X$. As explained in section~\ref{Sect11} below, this situation can be
generalized on the space of tensor fields of arbitrary type.

	It is worth writing also the (dual to~\eref{0.16}) implications (with
Lie derivatives)
    \begin{subequations}    \label{0.18}
    \begin{align}    \label{0.18a}
y'=0
&\ \Rightarrow\
Y' = \Bigl(\frac{\pd }{\pd U}\Bigr)' y
&\ \Rightarrow\
Y''-\tr C\cdot Y' + \det C\cdot Y = 0
&\ \Rightarrow\
Y_t = \frac{\pd }{\pd U} \e^{-Ct}y
\\		    \label{0.18b}
\varphi'=0
&\ \Rightarrow\
\Phi'= \varphi (\od U)'
&\ \Rightarrow\
\Phi''  - \tr C\cdot \Phi' + \det C\cdot \Phi = 0
&\ \Rightarrow\
\Phi_t = \varphi \e^{Ct} \od U
    \end{align}
    \end{subequations}
which describe the situation when the components $y$ and $\varphi$ are
invariant in the flow of $X$.

	An invariant of a linear vector field can easily be found by noticing
that the homothety operator
    \begin{equation}    \label{0.19}
P = u\frac{\pd }{\pd u} + v \frac{\pd }{\pd v}
    \end{equation}
commutes with $X$ and hence is an infinitesimal symmetry. One can calculate,
using the formalism from section~\ref{Sect0}, that the function
    \begin{equation}    \label{0.20}
I =
\frac{1}{2} \ln|W| -
\frac{c_1+c_2}{2} \int
\frac{\od p}{c_3+(c_1-c_4)p-c_2p^2}
    \end{equation}
with
\(
W=
    \begin{vmatrix}
u & u'\\
v & v'
    \end{vmatrix}
\)
being the Wronskian of $X$ and $p=\frac{v}{u}$is an invariant of $P$.


\section{Classification of the flows}
\label{Sect9}

		The classification of the flows and phase portraits of a linear
vector field on $\field[R]^2$ will be presented below by meas of the
eigenvalues $\lambda_1$ and $\lambda_2$ of the matrix $C$, which are generally
complex numbers. These are represented on the complex plane $\field[C]^2$ by
two points which form a segment with middle point
$\alpha=\frac{1}{2}(\lambda_1+\lambda_2)$. The following four cases are
possible:

\subparagraph{1. (Elliptical flow)}
	The eigenvalues $\lambda_1$ and $\lambda_2$ are complex conjugate, \ie
$\lambda_{1,2}=\alpha\pm\iu\beta$ for some $\alpha,\beta\in\field[R]$ with
$\beta\not=0$. An elliptic rotation is observed on the plane, affected by a
homothety for $\alpha\not=0$. At the origin, we have a stable (for $\alpha<0$)
or unstable (for $\alpha>0$) focus.

\subparagraph{2. (Hyperbolic flow)}
	The numbers $\lambda_1$ and $\lambda_2$ are different real ones,
$\lambda_{1,2}=\alpha\pm\beta$ for some $\alpha,\beta\in\field[R]$ with
$\beta\not=0$. Now a hyperbolic rotation can be seen on the plane which is
affected by a homothety for $\alpha\not=0$. At the origin we have a saddle
point for $\lambda_1\lambda_2<0$ or hyperbolic node for $\lambda_1\lambda_2>0$
which is stable if $\alpha<0$ or unstable if $\alpha>0$.

\subparagraph{3. (Parabolic flow)}
	The eigenvalues $\lambda_1$ and $\lambda_2$ are equal and
non-vanishing, $\lambda_{1}=\lambda_2=\alpha$ for some
$\alpha\in\field[R]\setminus\{0\}$ (and $\beta=0$). The flow is degenerate and
at the origin is observed a parabolic node which is stable if $\alpha<0$ or
unstable if $\alpha>0$.

\subparagraph{4. (Degenerate case)}
	The matrix $C$ is one-time degenerate, so that one of its eigenvalues
$\lambda_1$ or $\lambda_2$ vanishes
\[
\det C = 0 \quad \text{with} \quad
\lambda_1\lambda_2 = 0 \quad \lambda_1+\lambda_2\not=0 .
\]
Then there are numbers $a$ and $b$ such that one of them is non-zero and
\[
 (a,b)\cdot  \begin{pmatrix} c_1 & c_2 \\[1ex] c_3 & c_4 \end{pmatrix}
=0.
\]
Since~\eref{0.13} implies $au^\prime+bv^\prime=0$, the linear function
\[
I = au+bv
\]
is an invariant of the vector field $X$. The trajectories of $X$ are parallel
lines. Along them is observed an exponential movement of the points. For
instance, if $a=-1$ and $b=k$, the function
\[
 f=c_1u+c_2v
\]
is dragged via exponential law,
\[
f^\prime = kf
\quad\Leftrightarrow\quad
f_t = f\e^{kt}.
\]
The points situated on the line $f=\const$ move in an identical way and the
vector field $X$ is
\[
X = f \Bigl( \frac{\pd}{\pd x} + \frac{1}{k} \frac{\pd}{\pd v} \Bigr) .
\]

	We shall summarize the above considerations in the following
proposition.

\renewcommand{\thefigure}{\thesection.\arabic{figure}}

\begin{figure}[h!t!b!]
    \begin{center}
\hspace*{-3em}
\includegraphics[scale=0.72]{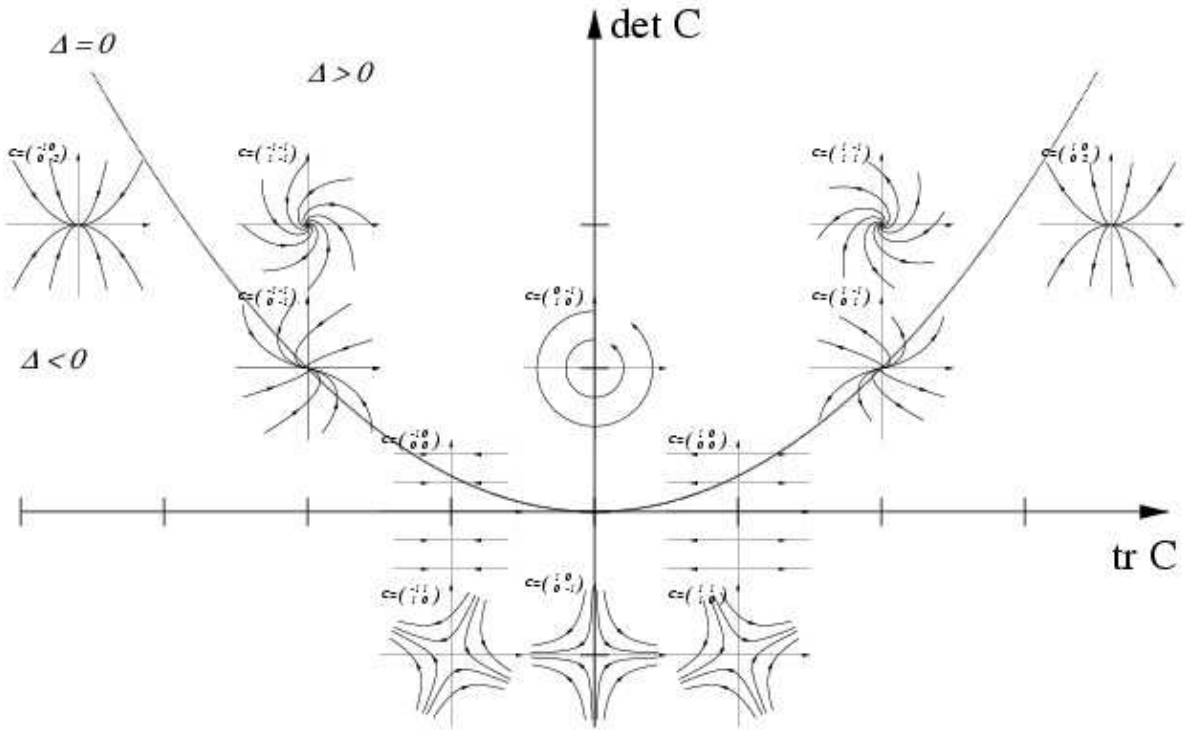}
\caption[Classification of the phase portraits]
	{\small Classification of the phase portraits on the plane $(x, y)$, with
	$x=\tr C$ and $y=\det C$, relative to the parabola
	$\Delta=x^2/4-y=0$.\label{fig9.2} }
    \end{center}
	\end{figure}

	On figure~\ref{fig9.2} are shown the phase portraits of a linear vector
field situated relative to the parabola
$\Delta=0$ with
\(
4 \Delta = (\tr C)^2 - 4\det C
        = (c_1-c_4)^2+4c_2c_3
        = (\lambda _{1}-\lambda_2)^2 .
\)
The focuses are inside the parabola
($\Delta<0$), the parabolic nodes ($\det C=0$) are on it ($\Delta=0$), and, at
last, the saddle points ($\det C<0$) and hyperbolic nodes ($\det C>0$) are
outside that curve ($\Delta>0$) (see also~\cite[p.~86]{Cantwell}, where a
similar picture is presented, but without some details).

    \begin{Prop}    \label{Prop8.2}
Let the system of inhomogeneous equations
	\begin{align*}
u^\prime &= c_1u+c_2v + b_1\\
v^\prime &= c_3u+c_4v + b_2 ,
	\end{align*}
or in a matrix form
\[
U^\prime = CU+B ,
\]
corresponds to a vector field $X$ (cf.~\eref{0.13}).
	If the matrix $C$ is nondegenerate, the flow of $X$ is described
similar to the one of a linear vector field
with the only difference that the singularity is
not at the coordinate origin but at a point $U_0$ such that
\[
CU_0+B = 0 .
\]
	The movement depends on is the rank of the matrix $(C|B)$, consisting
of the blocks $C$ and $B$, equal or not to the one of $C$.
    \end{Prop}

    \begin{Proof}
Take $U_0$ such that $CU_0+B=0$. Then the inhomogeneous system $U'=CU+B$ is
tantamount to the homogeneous system
\[
(U-U_0)' = C(U-U_0)
\]
relative to $(U-U_0)$ with a singularity at the point $U_0$.

	If the rank of $C$ is 2, the classification of
linear vector fields holds.
	When the rank of $C$ is 1, we can put without lost of generality
$c_1=kc_3$ and $c_2=kc_4$ for some number $k$. So that $X$  takes the form
\[
X =
  (c_3u+c_4v) \Bigl( k\frac{\pd}{\pd u} + \frac{\pd}{\pd v} \Bigr)
+ b_1\frac{\pd}{\pd u} + b_2\frac{\pd}{\pd v} .
\]
	The trajectories are exponential curves similar to the one of the
graph of the function $\field[R]\ni x\mapsto\e^x$.

	When $c_1=c_2=0$, we have only a uniform movement along  straight lines.
    \end{Proof}

\begin{figure}[h!tb]
    \begin{center}
\includegraphics[scale=1.3]{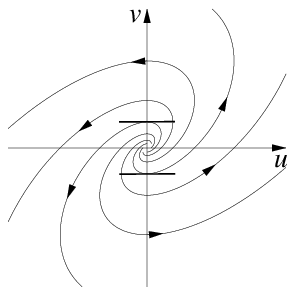}
\includegraphics[scale=1.3]{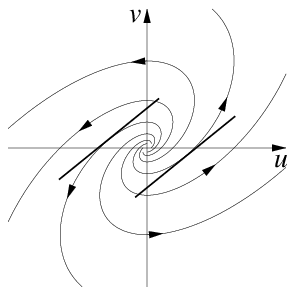}
\includegraphics[scale=1.3]{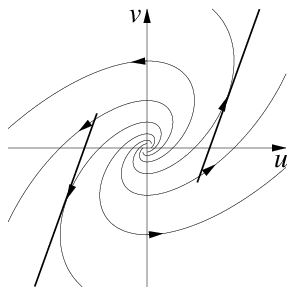}
\caption[Preservation of straight lines and parallelism]
	{\small Preservation of straight lines and
		parallelism.\label{fig9.3} }
    \end{center}
	\end{figure}

    \begin{Prop}    \label{Prop8.3}
The straight lines remain straight lines and their parallelism is preserved in
the flow of a linear vector field, as shown on figure~\ref{fig9.3}.
    \end{Prop}

    \begin{Proof}
Since the level lines of a linear function are straight lines, it
suffice to prove that a linear function is dragged by such a vector field into
a linear function. The flow is determined via the exponential flow
\[
U' = CU  \quad\Rightarrow\quad
U_t = \e^{Ct} U .
\]
Suppose $A=(a,b)$ is a matrix-raw  with $a,b\in\field[R]$. An arbitrary linear
function  $f$ has a representation  $f=AU$ and is dragged according to the law
$f_t=AU_t=A\e^{Ct}U$. Since the last function is linear for any fixed $t$, the
straight lines remain such and the parallelism is preserved; however, segments
of the lines can expand (for $\diver X>0$) or contract (for $\diver X<0$).
    \end{Proof}

\begin{figure}[h!tb]
	\begin{center}
 \includegraphics[scale=0.525]{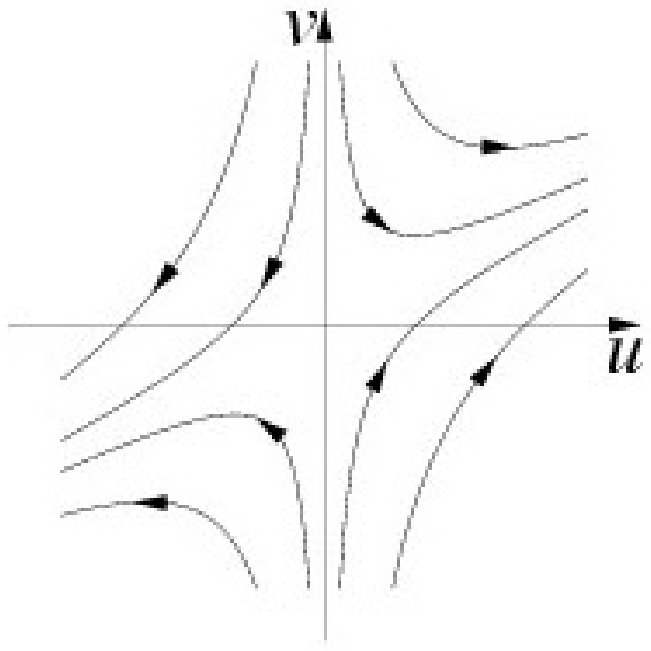}\hspace{2em}
 \includegraphics[scale=0.575]{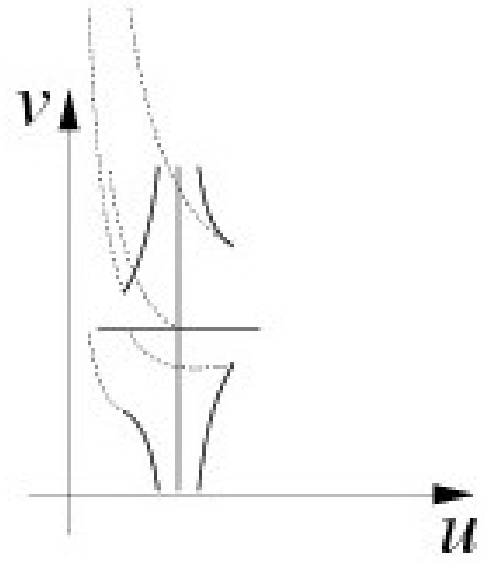}\hspace{2em}
 \includegraphics[scale=0.575]{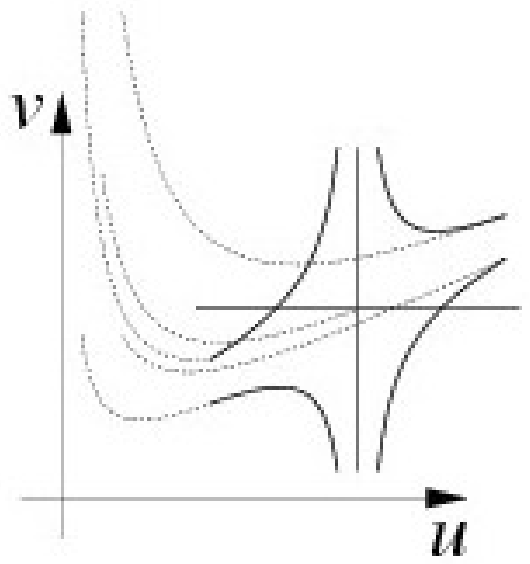}
\caption[Repetition of a phase portrait (saddle point) in moving frame]
	{\small Repetition of a phase portrait (saddle point) in moving frame.
	\label{fig9.4} }
	\end{center}
    \end{figure}

\begin{figure}[h!tb]
\begin{center}
\includegraphics[scale=1.6]{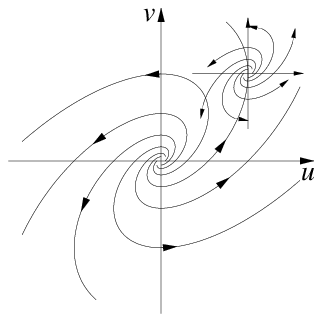}
\caption[Focus in a moving frame]
	{\small Focus in a moving frame.\label{fig9.5} }
	\end{center}
    \end{figure}

    \begin{Prop}    \label{Prop8.4}
When moving along the trajectories of a linear vector field together with a
coordinates system (moving frame), all points, moving along their own
trajectories, form in the moving frame the same phase portrait as the one in
which they are involved. This situation is illustrated on
figures~\ref{fig9.4} and~\ref{fig9.5}.
    \end{Prop}

    \begin{Proof}
The point $U$, where we are initially situated, is moving according to
$U_t=\e^{Ct}U$. At a neighboring point, we observe that it has radius\ndash
vector $U+\od U$ in the non-moving frame, but in the moving frame its
radius\ndash vector is $\od U$. From the law
\[
(U+\od U)_t = \e^{Ct}(U+\od U) ,
\]
we get $(dU)_t=\e^{Ct}\od U$, which proofs our assertion.~%
\footnote{~%
The same result follows also from $U'=CU$, \viz since $\od(U')=(\od U)'$, the
Lie derivative commutes with the exterior differentiation, then $(\od U)'=C\od
U$.%
}
The radius\ndash vectors $U$, relative to the non-moving frame, and $\od U$,
relative to the moving frame, are dragged by the flow in an identical way.
    \end{Proof}

    \begin{Rem}    \label{Rem9.1}
The Jacobi matrix  in a flow of a nonlinear vector field is non-constant and
can change when passing from one point to other.
	\begin{equation*}
    \begin{split}
u' &= u-v-u(u^2+v^2) \\
v' &= u+v-v(u^2+v^2)
    \end{split}
	\end{equation*}
with attractor $u^2+v^2=1$. To the five phase portraits correspond  Jacobi
matrices at five points when moving from the origin $(0,0)$ up along the $v$
axes.

	If an observer judges on a flow by taking into account only the Jacobi
matrix at a given point, the different observes obtain different results and
correspondingly they will have different opinions. In a case of a linear vector
field, the Jacobi matrix is constant and any local picture is an exact copy of
the global one.
\EndRem
    \end{Rem}


\section{Rotation groups in $\field[R]^3$}
\label{Sect10}

	The group of rotation of the space $\field[R]^3$ has three  basic
fundamental vector fields, say $X$, $Y$ and $Z$, which in a matrix notation can
be written as (in standard Cartesian coordinates $(u,v,w)$)
\[
(X,Y,Z)
=
\Bigl(\frac{\pd}{\pd u}, \frac{\pd}{\pd v}, \frac{\pd}{\pd w} \Bigr)
\cdot
\begin{pmatrix}
      0 & -w &  v \\
      w &  0 & -u \\
     -v &  u &  0
      \end{pmatrix} .
\]
As vector fields, the operators
    \begin{equation}    \label{10.1}
X =  w\frac{\pd}{\pd v} - v\frac{\pd}{\pd w} \quad
Y = -w\frac{\pd}{\pd u} + u\frac{\pd}{\pd w}\quad
Z =  v\frac{\pd}{\pd u} - u\frac{\pd}{\pd v}
    \end{equation}
are (linearly) dependent, \viz
\[
uX+vY+wZ = 0 ,
\]
and the function
\[
f = \frac{1}{2}(u^2+v^2+w^2)
\]
is their common invariant. The level surfaces of the function $f$ are
concentric spheres and are orbits of the rotation group. The table of
commutators of the operators $X$, $Y$ and $Z$ can easily be compted:
\renewcommand{\arraystretch}{0.2}
 	\begin{equation}    \label{10.2}
	\begin{tabular}{c|ccc}
$\nearrow$ & $X$ & $Y$ & $Z$\\
\hline\\
$X$ &0 & $Z$ & $-X$\\
&&&\\
$Y$ & $-Z$ & 0 & $Y$\\
&&&\\
$Z$ & $X$ & $-Y$ & 0  .
	\end{tabular}
 	\end{equation}

	A linear combination
\[
P = \omega_1X+\omega_2Y+\omega_3Z
\]
with constant coefficients $\omega_1,\omega_2,\omega_3\in\field[R]$ is also
a fundamental vector field of the rotation group, precisely of any
1\ndash parameter its subgroup. Besides the function $f$, an invariant of the
vector field
    \begin{equation}    \label{10.3}
P =
\begin{vmatrix}
        u & v & w \\[1ex]
	\omega_1 & \omega_2 & \omega_3\\[1ex]
        \frac{\pd}{\pd u} & \frac{\pd}{\pd v} & \frac{\pd}{\pd w}
\end{vmatrix}
    \end{equation}
is the linear function
\[
g = \omega_1u+\omega_2v+\omega_3w .
\]
Indeed, one can easily verify that $P(f)=P(g)=0$. This means that the
trajectories of $P$ are circles obtained as an intersection of the spheres
$f=\const$ with the planes $g=\const$; therefore $P$ is the rotation operator
around the axis $\omega_1:\omega_2:\omega_3$.

	The flow of $P$ drags the tensors and forces them to rotate around the
axis $\omega_1:\omega_2:\omega_3$. To illustrate that, we shall study the
dragging of $X$, $Y$ and $Z$ in the flow of $P$. A simple calculation of the
Lie derivatives relative to $P$, denoted via primes, gives (in a matrix form):
	\begin{align*}
(X,Y,Z)^\prime
&=
(X,Y,Z) \cdot
\begin{pmatrix}
        0 & -\omega_3 & \omega_2 \\[1ex]
        \omega_3 & 0 &-\omega_1 \\[1ex]
        -\omega_2 & \omega_1 & 0
        \end{pmatrix}
\\
(X,Y,Z)^{\prime\prime}
&=
(\omega_1,\omega_2,\omega_3) P - \omega (X,Y,Z)
\\
(X,Y,Z)^{\prime\prime\prime}
&=
  -\omega^2 (X,Y,Z)^\prime
\\
\omega^2
&:=
\omega_1^2+\omega_2^2+\omega_3^2.
	\end{align*}
It can be verified that anyone of the vector fields $X$, $Y$ and $Z$ satisfies
the equation
    \begin{gather}    \label{10.3-1}
S^{\prime\prime\prime} + \omega^2S^\prime = 0, \\
\intertext{whose general solution}    \label{10.3-2}
S_t = S + S^\prime\; \frac{\sin\omega t}{\omega}
        + S^{\prime\prime}\; \frac{1-\cos\omega t}{\omega^2}
    \end{gather}
reveals how these vector fields are dragged by the flow ot the field $P$.

	An interesting situation arises when the rotation is around the axis
1:1:1, \viz
\[
\omega_1=\omega_2=\omega_3=\frac{\sqrt{3}}{3} \quad
\omega=1 \quad
P=\frac{\sqrt{3}}{3}(X+Y+Z)
\]
\vspace{-5ex}
	\begin{multline*}
(X,Y,Z)_t
=
\frac{1}{3}(X+Y+Z) \cdot (1,1,1)
\\
 + \frac{2}{3}\Bigl[
(X,Y,Z) \cos t
+ (Z,X,Y) \cos\bigl( \frac{2\pi}{3}-t\bigr)
+ (Y,Z,X) \cos\bigl(\frac{2\pi}{3}+t\bigr)
\Bigr] .
	\end{multline*}
The three vector fields, occupying at $t=0$ the situation $(X,Y,Z)$, at
$t=\frac{2\pi}{3}$ move to $(Z,X,Y)$, and at $t=\frac{4\pi}{3}$ they move into
$(Y,Z,X)$. Therefore a cyclic permutation is in action.

	When the flows of $X$, $Y$ and $Z$ are dragged by the flow of $P$,
one should speak about a representation of the Lie algebra of the rotation
group or of the adjoint representation of the rotation group.

	Let us generalize the above setting~\cite{Rahula&Tseluiko}. Consider
the following quadratic forms on $\field[R]^3$ coordinated by $(u,v,w)$:
	\begin{equation}    \label{10.5-1}
    \begin{split}
f &= \frac{1}{2}(a_{11}u^2+a_{22}v^2+a_{33}w^2 +2a_{12}uv+2a_{13}uw+2a_{23}vw)
\\
\tilde{f} &= \frac{1}{2}(\bar{a}_{11}u^2+\bar{a}_{22}v^2
             + \bar{a}_{33}w^2+2\bar{a}_{12}uv+2\bar{a}_{13}uw+2\bar{a}_{23}vw).
    \end{split}
	\end{equation}
The matrices of their coefficients
\[
A= \begin{pmatrix}
        a_{11}&a_{12}&a_{13}\\[1ex]
        a_{12}&a_{22}&a_{23}\\[1ex]
        a_{13}&a_{23}&a_{33}\\[1ex]
        \end{pmatrix}
\qquad
\bar{A}= \begin{pmatrix}
        \bar{a}_{11}&\bar{a}_{12}&\bar{a}_{13}\\[1ex]
        \bar{a}_{12}&\bar{a}_{22}&\bar{a}_{23}\\[1ex]
        \bar{a}_{13}&\bar{a}_{23}&\bar{a}_{33}\\[1ex]
        \end{pmatrix}
\]
are required to be symmetric and mutually inverse, $A^\top=A$,
$\bar{A}^\top=\bar{A}$ and $A\cdot\bar{A}=\openone$, so that $(\det
A)(\det\bar{A})=1$.

	The level surfaces of $f$ and $\bar{f}$ are central quadrics, \viz
ellipsoids, when $f$ and $\bar{f}$ have constant signs, or hyperboloids, when
$f$ and $\bar{f}$ can change signs. The equation $\bar{f}=\const$ is called
tangent for the surface $f=\const$.

	Let us introduce the following shortcuts for the partial derivatives of
$f$:
\[
 f_i=a_{i1}u+a_{i2}v+a_{i3}w,\quad i=1,2,3.
\]

    \begin{Prop}    \label{Prop10.1}
The vector fields $X_1$, $X_2$ and $X_3$ given by
    \begin{equation}    \label{10.4}
(X_1,X_2,X_3)
=
\Bigl( \frac{\pd}{\pd u}, \frac{\pd}{\pd v}, \frac{\pd}{\pd w} \Bigr)
\cdot
\begin{pmatrix}
      0    &-f_3 &f_2\\[1ex]
      f_3  &0    &-f_1\\[1ex]
      -f_2 &f_1  &0
      \end{pmatrix}
    \end{equation}
are fundamental vector fields of a 3-dimensional Lie group $G$ and
 	\begin{equation}    \label{10.5}
	\begin{tabular}{c|ccc}
$\nearrow$ & $X_1$ & $X_2$ & $X_3$ \\
\hline\\
$X_1$ & $0$ &$\sum_{k=1}^3a_{3k}X_k$  &       $-\sum_{k=1}^3a_{2k}X_k$ \\[2ex]
$X_2$ &      $-\sum_{k=1}^3a_{3k}X_k$ & $0$ & $\sum_{k=1}^3a_{1k}X_k$ \\[2ex]
$X_3$ &      $\sum_{k=1}^3a_{2k}X_k$  &       $-\sum_{k=1}^3a_{1k}X_k$  & $0$
	\end{tabular}
 	\end{equation}
is the table of their commutators. The fundamental vector fields given
by~\eref{10.4} are (linearly) dependent,
\[
f_1X_1+f_2X_2+f_3X_3 = 0,
\]
and the family of quadrics $f=\const$ are orbits of the group $G$.
    \end{Prop}

    \begin{Proof}
The commutation table is obtained via direct calculations like
	\begin{align*}
[X_1,X_2]
&=
\Bigl[ f_3\frac{\pd}{\pd v} - f_2 \frac{\pd}{\pd w}
,
       f_1\frac{\pd}{\pd w} - f_3\frac{\pd}{\pd u}
\Bigr]
\\
&= (a_{33}f_2-a_{23}f_3)\frac{\pd}{\pd u}+
   (a_{13}f_3-a_{33}f_1)\frac{\pd}{\pd v}+
   (a_{23}f_1-a_{13}f_2)\frac{\pd}{\pd w}
\\
&= a_{31}X_1+a_{32}X_2+a_{33}X_3 .
	\end{align*}

According the the third Lie theorem
(see~\cite[p.~283 of the English text]{Mathenedia-3}
or~\cite[p.~102]{Chebotarev-1940}), the coefficients of the quadratic form $f$
define the structure constants $c_{jk}^{i}$ of the group $G$, \viz
	\begin{alignat*}{2}
a_{11}=c_{23}^1  &&\qquad& a_{12}=c_{23}^2=c_{31}^1 \\
a_{22}=c_{31}^2  &&& a_{23}=c_{12}^2=c_{31}^3 \\
a_{33}=c_{12}^3  &&& a_{13}=c_{12}^1=c_{23}^3,
	\end{alignat*}
and the fundamental vector fields of $G$ in $\field[R]^3$, coordinated by
$(u,v,w)$, are the vector fields $X_1$, $X_2$ and $X_3$. The orbits of $G$ are
the quadrics $f=\const$ as $X_i(f)=0$, $i=1,2,3$, and the rank of the matrix
\[
\begin{pmatrix}
        0&-f_3&f_2\\[1ex]
        f_3&0&-f_1\\[1ex]
        -f_2&f_1&0
        \end{pmatrix}
\]
equals 2; in the opposite case, we should have $f_1=f_2=f_3=0$ which will
mean that the matrix $A$ is degenerate.
    \end{Proof}

    \begin{Prop}    \label{Prop10.2}
The vector field
    \begin{equation}    \label{10.6}
P=p_1X_1+p_2X_2+p_3X_3,
    \end{equation}
which is a linear combination of $X_1$, $X_2$ and $X_3$ with constant
coefficients $p_1,p_2,p_3\in\field[R]$ (and hence is a fundamental field of
$G$), of $G$ has the invariant
\[
g=p_1u+p_2v+p_3w
\]
besides the function $f$.

	The trajectories of $P$ are on the intersection of the quadrics
$f=\const$ with the planes $g=\const$ and, depending on the sign of the
determinant
    \begin{equation}    \label{10.7}
\varepsilon
:=\begin{vmatrix}
            a_{11} &a_{12} &a_{13} &p_1\\[1ex]
            a_{12} &a_{22} &a_{23} &p_2\\[1ex]
            a_{13} &a_{23} &a_{33} &p_3\\[0.5ex]
            p_1    &p_2    &p_3    &0
	\end{vmatrix},
    \end{equation}
define:
	\begin{itemize}
\item
	a family of ellipses, including isolated points, for $\varepsilon<0$.
\item
	a family of hyperbolas, including pairs of intersecting straight
lines, for $\varepsilon>0$.
\item
	a family of parabolas, including pairs of parallel or coinciding
straight lines, for $\varepsilon=0$.
	\end{itemize}
    \end{Prop}

    \begin{Rem}    \label{Rem10.0-1}
	The possible flows on the quadrics $f=\const$ are depicted on
figures~\ref{fig10.1} and~\ref{fig10.4}.

	\begin{figure}[h!tb]
    \begin{center}
	\begin{tabular}[h!t]{cc}
\includegraphics[scale=1.1]{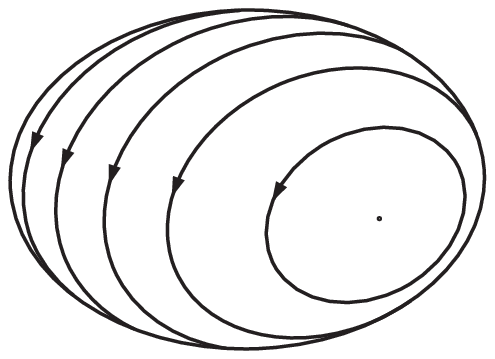}
&  \hspace{2.71828em}
\includegraphics[scale=0.7]{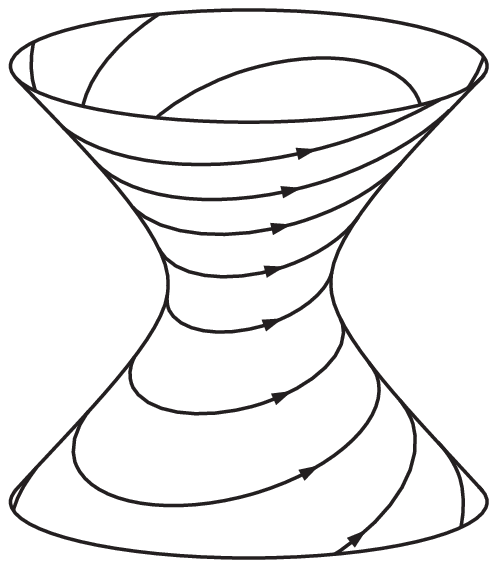}
	\end{tabular}
\caption[Elliptical flow on an ellipsoid and on a hyperboloid]
	{\small Elliptical flow on an ellipsoid (left) and on a hyperboloid
	(right).\label{fig10.1} }
    \end{center}
	\end{figure}
\vspace{5ex}
	\begin{figure}[h!tb]
    \begin{center}
	\begin{tabular}[h!t]{cc}
\includegraphics[scale=0.8]{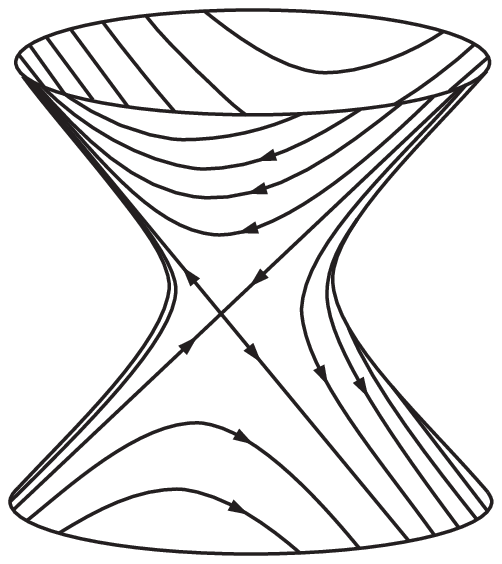}
&  \hspace{2.71828em}
\includegraphics[scale=0.8]{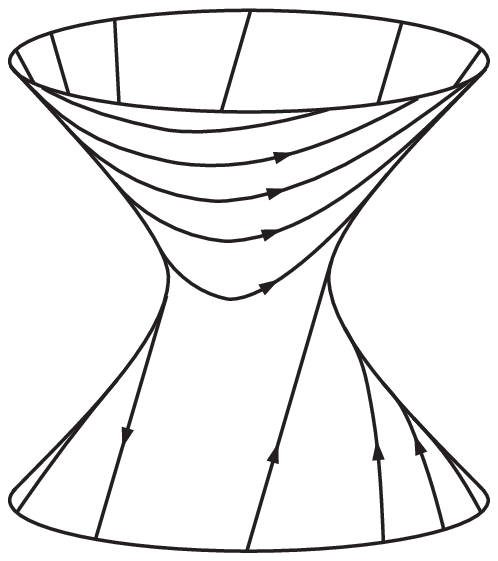}
	\end{tabular}
\caption[Hyperbolic and parabolic flow on a hyperboloid]
	{\small Hyperbolic (left) and parabolic (right) flow on a
	hyperboloid.\label{fig10.4} }
    \end{center}
	\end{figure}

    \end{Rem}

    \begin{Proof}
Representing the field $P$ as
\[
P = \begin{pmatrix}
	f_1 & f_2 & f_3 \\[1ex]
	u   & v   & w   \\[1ex]
        \frac{\pd}{\pd u} & \frac{\pd}{\pd v} & \frac{\pd}{\pd w}
   \end{pmatrix} ,
\]
we immediately get $P(f)=P(g)=0$, \ie $f$ and $g$ are invariants of $P$.
Therefore the trajectories of $P$ lie on the intersection of the quadrics
$f=\const$ with the planes $g=\const$.

	The type of the trajectories is determined by the intersection of the
plane $g=0$ with the cone $f=0$, \ie by the system $f=g=0$. If this system has
only the zero solution $u=v=0$, the cone $f=0$ can be either imaginary and the
surfaces $f=\const$ are ellipsoids or real and the surfaces $f=\const$ are
hyperbolloids. In the both cases a family of ellipsoids is observed on the
planes $g=\const$, as shown on figure~\ref{fig10.1}.

	If the system defines a pair of different lines, the plane $g=0$
intersects the real cone $f=0$ along two generants and, on the plane $g=0$, one
can observe a family of hyperbolas (see figure~\ref{fig10.4}, left). If the
systems defines two identical lines, the plane $g=0$ is tangent to the real
tangential cone $\tilde{f}=0$ and on the planes $g=\const$ one can see a family
of parabolas (see figure~\ref{fig10.4}, right).
	This is equivalent that the vector $p_1:p_2:p_3$, from the origin, to be
situated inside the tangent cone $f=0$, outside this cone, or on the cone.

	Analytically the situation is as follows. We suppose that $p_3\not=0$
without loss of generality. From $g=0$, we get $p_3w=-p_1u-p_2v$, so that the
substitution into $p_3^2f$ results in
	\begin{multline*}
(a_{11}p_3^2-2a_{13}p_1p_2+a_{33}p_1^2)u^2
  + 2(a_{33}p_1p_2-a_{13}p_2p_3-a_{23}p_1p_3+a_{12}p_3^2)uv
\\
  + (a_{22}p_3^2-2a_{13}p_2p_3+a_{33}p_2^2)v^2.
	\end{multline*}
and the corresponding equadratic equation defines the fraction $u/v$.

	In this way one gets the relation $u:v:w$ and hence the intersection
of the plane $g=0$ with the cone $f=0$. All depends on the discriminant which,
up to the positive coefficient $p_3^2$, coincides with the quantity
$\varepsilon$, \viz
	\begin{multline*}
(a_{33}p_1p_2-a_{13}p_2p_3-a_{23}p_1p_3+a_{12}p_3^2)^2
\\
- (a_{11}p_3^2-2a_{13}p_1p_2+a_{33}p_1^2)(a_{22}p_3^2-2a_{13}p_2p_3+a_{33}p_2^2)
\\
=p_3^2\Big(
  -\begin{vmatrix}a_{11} & a_{12} \\[0,5ex] a_{12} & a_{22} \end{vmatrix} p_3^2
  -\begin{vmatrix}a_{22} & a_{23} \\[0,5ex] a_{23} & a_{33} \end{vmatrix} p_1^2
  -\begin{vmatrix}a_{11} & a_{13} \\[0,5ex] a_{13} & a_{33} \end{vmatrix} p_2^2
 +2\begin{vmatrix}a_{12} & a_{13} \\[0,5ex] a_{23} & a_{33} \end{vmatrix} p_1p_2
\\[0,5ex]
 +2\begin{vmatrix} a_{13} & a_{12} \\[0,5ex] a_{23} & a_{22} \end{vmatrix} p_1p_3
 +2\begin{vmatrix} a_{11} & a_{12} \\[0,5ex] a_{13} & a_{23} \end{vmatrix} p_2p_3
\Big)
=p_3^2\varepsilon.
	\end{multline*}
	Consequently, for $\varepsilon<0$ (for $\varepsilon>0$) the flow of $P$
is elliptic (hyperbolic) and for $\varepsilon=0$ it is parabolic.
    \end{Proof}

    \begin{Rem}    \label{Rem10.1}
Generally, we have a class of Lie groups $G$ each of which has its own
 table of commutators and structure constants. In particular, this may be the
classical rotation group with the coefficients
\[
a_{11} = a_{22} = a_{33}=1 \quad a_{12} = a_{13} = a_{23} = 0
\]
or the function
\[
f = \frac{1}{2}(u^2+v^2+w^2) .
\]

	A more interesting group is determined via the table of
commutators
\renewcommand{\arraystretch}{0.2}
 	\begin{equation}    \label{7.6}
  \begin{tabular}{c|r@{\hspace{-0.1em}}c@{\hspace{-0.2em}}ccc}
$\nearrow$ & $X_0$&$\vdots$ & $X_1$ &$X_2$ &$X_3$ \\[0.4ex]
\hline\\\vspace*{0.6ex}
$X_0$	   & 0&$\vdots$ &0&0&0\\
\dots &\dots &\dots &\dots &\dots &\dots  \\
$X_1$	   & 0&$\vdots$ &0&$2X_3$&$2X_2$\\[0.4ex]
$X_2$	   & 0&$\vdots$ &$-2X_3$&0&$-2X_1$\\[0.4ex]
$X_3$	   & 0&$\vdots$ &$-2X_2$&$2X_1$&0 ,
  \end{tabular}
 	\end{equation}
in which the value of the commutator $[X_i,X_j]:=X_i\circ X_j-X_j\circ X_i$,
$i,j=0,1,2,3$, is situated at the intersection between the $i^{\text{th}}$ row
and $j^{\text{th}}$ column,
and, more precisely by its part which is related to
the factor-group of the group of centriaffine transformations relative to the
homothetic subgroup, \ie to the group of equiaffine transformations of the
plane. This group is distinguished by the coefficients
\[
a_{33}=-a_{22}-a_{11}=2 \quad  a_{12}=a_{13}=a_{23}=0 ,
\]
or via the function
\[
f = -(u^2+v^2-w^2) .
\]
Of course, this is the group of pseudo-Euclidean transformations of
$\field[R]^3$, \ie the Lorentz group in three dimensions.
\EndRem
    \end{Rem}

    \begin{Rem}    \label{Rem10.2}
One can observe a circular movements on a sphere which themselves are involved in
a rotational movement around some axis; \eg the dragging of a field $X$ by the
flow of a field $P$. The notion of a ``transformation of a transformation'' is
not new; see, e.g.,~\cite[p.~26]{Chebotarev-1940}. Here we can talk of
``movements of a movement'', ``movements of the movements of a movement'', etc.,
\ie of movements of higher orders., which seems to be a good item for a future
research. For instance, one can imagine that a given flow undergoes a
transformation or that it is dragged by the flow of other field, that this
process undergoes a transformation or is dragged by the flow of a third field,
\etc
\EndRem
    \end{Rem}

	One can easily imagine an elliptic flow on an ellipsoid, which is
dragged by other elliptic flow. A hyperboloid can be intersected by a family of
parallel planes, the result being a family of ellipses, or hyperbolas, or
parabolas. Consequently on a hyperboloid can be observed an elliptic,
hyperbolic and parabolic flows. However, it is difficult to be imagined the
change on a hyperboloid, when an elliptic flow is continuously deformed and
that at some moment it transforms into a parabolic flow and, then, into a
hyperbolic; similarly, a stable hyperbolic flow can change into  a parabolic and
then into a stable elliptic flow. Such bifurcations in a (single or two sided)
hyperboloid are admissible.

	Consider now how the (flow of the) vector field $P$ influences the
operators $X_1$, $X_2$ and $X_3$.

    \begin{Prop}    \label{Prop10.3}
The Lie derivatives  $X_i^\prime = [P,X_i]$, $i=1, 2, 3$, of the fundamental
vector fields $X_1$, $X_2$ and $X_3$ relative to the field $P$ can be expressed
through the same fields via the equation
    \begin{equation}    \label{10.8}
(X_1,X_2,X_3)^\prime
=
(X_1,X_2,X_3) \cdot B ,
    \end{equation}
with $B$ being the following product matrix
\[
B
=
\begin{pmatrix}	a_{11} & a_{12} & a_{13}\\[1ex]
		a_{12} & a_{22} & a_{23}\\[1ex]
		a_{13} & a_{23} & a_{33} \end{pmatrix}
\cdot
\begin{pmatrix}	0 	& -p_3 & p_2\\[1ex]
		p_3  	& 0    &-p_1\\[1ex]
		-p_2 	& p_1  & 0	\end{pmatrix}
\]
which is such that
\[
B^3 = \varepsilon B ,
\]
where $\varepsilon$ is given by~\eref{10.7}. Besides, any one of the
fundamental vector fields $X_1$, $X_2$ and $X_3$ is a solution of the
differential equation (cf.~\eref{10.3-1})

\[
S^{\prime\prime\prime} - \varepsilon S^\prime = 0 ,
\]
which, depending on the sign of $\varepsilon$, admits the following solutions
$S_t$:
	\begin{itemize}
     \item[1)] If $\varepsilon<0$ and $\lambda=\sqrt{-\varepsilon}$, then
(cf.~\eref{10.3-2})
\[
S_t = S + S^\prime	\frac{\sin\lambda t}{\lambda}
      + S^{\prime\prime}\frac{1-\cos\lambda t}{\lambda^2} .
\]

     \item[2)] If $\varepsilon>0$ and $\lambda=\sqrt{\varepsilon}$, then
\[
S_t = S + S^\prime	\frac{\sinh \lambda t}{\lambda}
      + S^{\prime\prime}\frac{1-\cosh\lambda t}{\lambda^2} .
\]

     \item[3)] If $\varepsilon=0$, then
\[
S_t = S + S^\prime t + S^{\prime\prime} \frac{t^2}{2} .
\]
	\end{itemize}
    \end{Prop}

    \begin{Proof}
The Lie derivative $X'_1$ is calculated by using the commutation
table~\eref{10.5}:
	\begin{multline*}
X_1^\prime
:=[P,X_1]
= [p_1X_1+p_2X_2+p_3X_3,X_1]
= p_2[X_2X_1]+p_3[X_3X_1]
\\
= (a_{12}p_3-a_{13}p_1)X_1+(a_{22}p_3-a_{23}p_1)X_2+(a_{23}p_3-a_{33}p_1)X_3
\\
= (X_1,X_2,X_3)\cdot \begin{pmatrix} a_{11}&a_{12}&a_{13}\\[1ex]
				     a_{12}&a_{22}&a_{23}\\[1ex]
				     a_{13}&a_{23}&a_{33}\end{pmatrix}
\cdot \begin{pmatrix} 0 \\[1ex] p_3\\[1ex] -p_2 \end{pmatrix} .
	\end{multline*}
The Lie derivatives $X'_2$ and $X'_3$ can be calculated similarly.

	Since the Hamilton-Cayley formula for the matrix $B$ reads
\[ \sigma_0 B^3 + \sigma_1 B^2 + \sigma_2 B + \sigma_3\openone = 0 ,
\]
where $\sigma_0=1$, $\sigma_1=\tr B=0$ and $\sigma_3=\det B=0$, we have
\[
B^3 + \sigma_2 B = 0 .
\]
Calculating $\sigma_2$,
	\begin{multline*}
\sigma_2
=
\begin{vmatrix} b_{11} & b_{12} \\[1ex] b_{12} & b_{22} \end{vmatrix} +
\begin{vmatrix} b_{22} & b_{23} \\[1ex] b_{23} & b_{33} \end{vmatrix} +
\begin{vmatrix} b_{11} & b_{13} \\[1ex] b_{13} & b_{33} \end{vmatrix}
\\
=
  (a_{12}p_3-a_{13}p_2) (a_{23}p_1-a_{12}p_3)
- (a_{23}p_1-a_{12}p_3) (a_{13}p_2-a_{23}p_1)
\\
+ (a_{12}p_3-a_{13}p_2) (a_{13}p_2-a_{23}p_1)
- (a_{13}p_1-a_{11}p_3)(a_{22}p_3-a_{23}p_2)
\\
- (a_{11}p_2-a_{12}p_1) (a_{23}p_3-a_{33}p_2)
- (a_{12}p_2-a_{22}p_1)(a_{33}p_1-a_{13}p_3) ,
	\end{multline*}
we find $\sigma_2=-\varepsilon$.

	The establishment of the equation
$S^{\prime\prime\prime}-\varepsilon S^\prime = 0$
for the fields $X_1$, $X_2$ and $X_3$ is trivial. A simple verification reveals
that the written functions $S_t$ are its solutions. Besides, the third solution
is obtained from any one of the preceding two in the limit $\lambda\to0$ due to
\[
\lim_{\lambda\to 0}\frac{\sin\lambda t}{\lambda}
=
 \lim_{\lambda\to 0}\frac{\sinh\lambda t}{\lambda}=t
\qquad
\lim_{\lambda\to 0}\frac{1-\cos\lambda t}{\lambda^2}=
  \lim_{\lambda\to 0}\frac{1-\cosh\lambda t}{\lambda^2}=\frac{t^2}{2} \, .
\]
    \end{Proof}

    \begin{Rem}    \label{Rem10.3}
The behavior of the vector fields $X_1$, $X_2$ and $X_3$ in the flow of $P$ is
consistent with that flow in a sense that, if the flow is elliptic, hyperbolic
or parabolic (see the three cases in proposition~\ref{Prop10.2}), then the
behavior of anyone of these fields is respectively elliptic, hyperbolic or
parabolic (see cases~1,~2 and~3 of proposition~\ref{Prop10.3}). A
bifurcation can occur if the sign of $\varepsilon$ can change or, in
geometrical language, if the normal to the plane $g=0$ with direction
$p_1:p_2:p_3$ is moved from the inner part of the tangent cone $\bar{f}=0$ to
the outside domain, or \emph{vice versa} if it moves from the outside part of the
cone to its interior.
\EndRem
    \end{Rem}

	A similar bifurcation can occur with the flow of a linear vector field
$X$ on the plane, \ie a transition form elliptic to hyperbolic behavior or
\emph{vice versa}, if the quantity $\Delta$ changes its sign (see cases 1, 2 or
3 of proposition~\ref{Prop10.3}). The equality $\Delta=0$ in the space
$\field[R]^3$ coordinated by $b_1$, $b_2$ and $b_3$ should be thought as an
equation of the cone
\[
b_1^2 + b_2^2 - b_3^2 = 0
\]
and the bifurcation is due to the change     of the direction $b_1:b_2:b_3$
from inside to outside the cone or from outside to inside of it.


\section
[Tensor representations of the group \protect{$\GL(2,\field[R])$} and the
algebra \protect{$\gl(2,\field[R])$}]
{Tensor representations of the group $\GL(2,\field[R])$ and \\ the algebra
 $\gl(2,\field[R])$}
\label{Sect11}

Let on a $2$-dimensional manifold be given a frame $R$ and its dual coframe
$\Theta$. In them a tensor $S$ of type $(p,q)\in\field[N]^2$ is represented
according
	\begin{equation}    \label{4.11}
S
=R_{i_{1}}\otimes\dots\otimes R_{i_{p}}
 \otimes s^{i_{1}\dots i_{p}}_{j_{1}\dots j_{q}} \otimes
 \theta ^{j_{1}}\otimes \dots \otimes \theta ^{j_{q}} .
	\end{equation}
whith
\(
s_{j_1\dots j_q}^{i_1\dots i_p}
=
S(\theta^{i_1},\dots,\theta^{i_p}; R_{j_1},\dots,R_{j_q})
\)
being the components of $S$
relative to the tensor frame induced by $\{R_i\}$ and $\{\theta^i\}$.
	A regular $2\times 2$ matrix
$A\in\GL(2,\field[R])$ transforms S into a tensor $\tilde{S}$ with local
components
	\begin{equation}    \label{11.1}
\tilde{s}_{j_1...j_q}^{i_1...i_p}
=
A_{k_1}^{i_1}...A_{k_p}^{i_p}s_{l_1...l_q}^{k_1...k_p}
\bar{A}_{j_1}^{l_1}...\bar{A}_{j_p}^{l_q},
	\end{equation}
where $A_j^i$ and $\bar{A}_j^i$ are the elements of $A$ and its inverse
matrix $A^{-1}$, respectively. It is said that the group $\GL(2,\field[R])$ is
represented or that it acts via~\eref{11.1} on the space of tensors of type
$(p,q)$.

	An arbitrary matrix $C$ of type $2\times2$ is an element of the Lie
algebra, which is isomorphic with the Lie algebra $\gl(2,\field[R])$ of the
group $\GL(2,\field[R])$. The exponential $\e^{Ct}$, $t\in\field[R]$, defines a
1\ndash parameter subgroup of $\GL(2,\field[R])$ and equation~\eref{11.1}
with $A=\e^{Ct}$ defines a 1\ndash parameter family of tensors $S_t$. In this
way, in the space of tensors of type $(p,q)$, is defined the flow of some
vector field $\bar{X}$ whose local components can be found by putting
$A=\e^{Ct}$ in~\eref{11.1}, differentiating the result with respect to $t$
and putting $t=0$. To the vector field $\bar{X}$
corresponds a matrix $\bar{C}$, similarly to the correspondence $X\to C$; the
structure of $\bar{C}$ can be obtained by ordering the components of the tensor
$S$. We say that the correspondence $C\to\bar{C}$ or $X\to\bar{X}$ determines a
representation of the Lie algebra $\gl(2,\field[R])$ on the space of tensors
$S$.

	It is important to find the eigenvalues of $\bar{C}$ provided the ones
of $C$ are known. For instance, with their help one can restore the
exponential $\e^{Ct}$.
The flow of $\bar{X}$ is determined directly by~\eref{11.1} after the
substitution $A=\e^{Ct}$.

    \begin{Prop}    \label{Prop11.1}
Let the group $\GL(2,\field[R])$ acts on the space of tensors of
type $(p,q)$ on a 2\ndash dimensional manifold. If $\lambda_1$ and $\lambda_2$
are the eigenvalues of $C$, then the ones of $\bar{C}$ are
	\begin{equation}    \label{11.2}
\lambda_{i_1} + \dots + \lambda_{i_p} - \lambda_{j_1} - \dots - \lambda_{j_q} ,
	\end{equation}
where all indices take the values~1 and~2. The number of these eigenvalues
equals the number of the components of a tensor of type $(p,q)$, \ie to
$2^{p+q}$  in two dimensions.
    \end{Prop}

    \begin{Proof}
Chose a frame $R$ such that $C$ will take in it the diagonal form
$C=\diag(\lambda_1,\lambda_2)$ with $\lambda_1,\lambda_2\in\field[C]$
being the eigenvalues of $C$. Then
$\e^{Ct}=\diag(\e^{\lambda_1t},\e^{\lambda_2t})$, so that~\eref{11.1} with
$A=\e^{Ct}$ reads
	\begin{equation}
\e^{\lambda_{i_1}t} \dots \e^{\lambda_{i_p}t}
  s_{j_1\dots j_q}^{i_1...i_p}
\e^{-\lambda_{j_1}t} \dots \e^{-\lambda_{j_q}t}
=
  s_{j_1...j_q}^{i_1...i_p}
\e^{( \lambda_{i_1}+ \dots +\lambda_{i_p}
     -\lambda_{j_1}- \dots-\lambda_{j_q})t} .
	\end{equation}
Differentiating this expression with respect to $t$ and setting $t=0$, we get
the components of $S$ multiplied by the numbers~\eref{11.2}. Thus we get a
diagonal matrix with diagonal elements~\eref{11.2}, which are the
eigenvalues of $\bar{C}$.
    \end{Proof}

\begin{figure}
    \begin{center}
\includegraphics[scale=0.525]{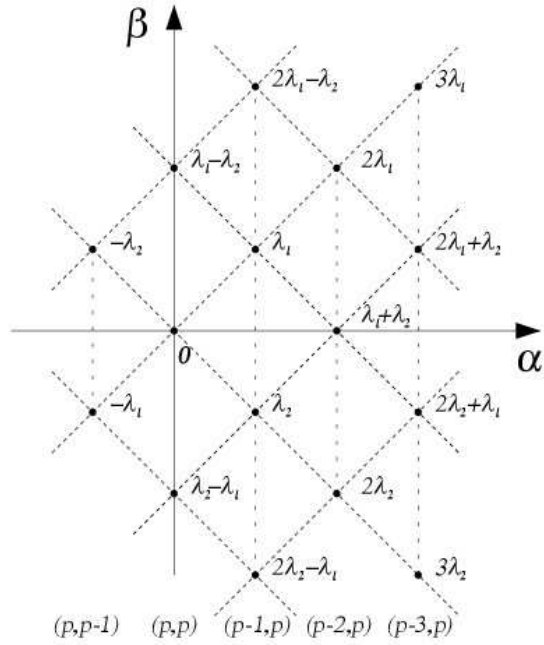}
\caption[Distribution of eigenvalues on the complex plane]
	{\small Distribution of the numbers~\protect{\eref{11.2}} on the complex
	plane for $\lambda_{1,2}=\alpha\pm\iu\beta$.\label{fig11.1} }
    \end{center}
	\end{figure}

	One can restore the linear ordinary differential equation for
the components of the tensor field $S$ (which is invariant in the given flow)
via the quantities~\eref{11.2}. The coefficients are symmetric polynomials
from~\eref{11.2}. Similarly can be restored the differential equation for
the tensor field $S$ with invariant components, but via the quantities with
inverse signs according to the duality principle, \viz
	\begin{equation}    \label{11.3}
  \lambda_{j_1} + \dots + \lambda_{j_q}
- \lambda_{i_1} - \dots - \lambda_{i_p}.
	\end{equation}
The numbers~\eref{11.2} and~\eref{11.3} are situated on the complex plane
on the knots of some lattice determined by the initial values $\lambda_1$ and
$\lambda_2$. On figure~\ref{fig11.1}
is shown the position of the
numbers~\eref{11.3} on the base of the complex conjugate eigenvalues
$\lambda_{1,2}=\alpha\pm\iu\beta$.
%
%

	When $p$ and $q$ are fixed, the numbers~\eref{11.3} have identical
real parts, equal to $(q-p)\alpha$, and, consequently, are situated on a
straight line parallel to the imaginary axis.

	Let us consider three particular cases of tensor fields dragged by the
flow of a linear vector field $X$:

	To illustrate the general situation, below we shall reproduce three
propositions from~\cite[sec.~2.9]{Rahula-VectorFields}, where the corresponding
proofs can be found.

    \begin{Exmp}    \label{Exmp11.2}
	For a tensor field of type $(1,1)$ (an affinor)
    \begin{equation}    \label{11.4}
S= s_j^iR_i\otimes\theta^j ,
    \end{equation}
where $s=[s_j^i]$ is $2\times2$ matrix, we have the eigenvalues
	\begin{equation}    \label{11.5}
\lambda_1-\lambda_2 , \quad
\lambda_2-\lambda_1 , \quad
0 \text{ (a double eigenvalue)} .
	\end{equation}
For it are valid the following
implications
	\begin{alignat}{2}    \label{11.11}
& s^\prime=0 &\ \Rightarrow\ &
	S^{\prime\prime\prime} = \Delta\cdot S^\prime
\\     \label{11.12}
& S^\prime=0 &\ \Rightarrow\ &
s^{\prime\prime\prime} = \Delta\cdot s^\prime ,
	\end{alignat}
with $\Delta$ defined by
	\begin{equation}    \label{5.1}
	\begin{split}
4 \Delta &= (\tr C)^2 - 4\det C
        = (c_1-c_4)^2+4c_2c_3
        = (\lambda _{1}-\lambda_2)^2 .
	\end{split}
	\end{equation}
Besides, in a case of~\eref{11.12}, the matrix
$\bar{C}$ is
	\begin{equation}    \label{11.13}
\bar{C} =
\begin{pmatrix}
0   & -c_3    & c_2  & 0    \\[1ex]
-c_2& c_1-c_4 &     0& c_3 \\[1ex]
c_3 & 0       & c_4-c_1 &-c_2 \\[1ex]
0   &c_3      &-c_2  & 0
\end{pmatrix}
	\end{equation}
and the vector field $\bar{X}$ is given by
\[
\bar{X} =
\Bigl(
\frac{\pd }{\pd s_1^1}, \frac{\pd }{\pd s_2^1} ,
\frac{\pd }{\pd s_1^2}, \frac{\pd }{\pd s_2^2}
\Bigr)
\cdot \bar{C} \cdot (s_1^1,s_2^1,s_1^2,s_2^2)^\top .
\]
    \end{Exmp}


    \begin{Exmp}    \label{Exmp11.3}
	For tensor fields of types $(0,2)$ and $(0,3)$, \ie for quadratic and
cubic forms
    \begin{subequations}	\label{11.6}
	\begin{align}    \label{11.6a}
G &= g_{ij} \theta^i \otimes \theta^j
\\			 \label{11.6b}
H &= h_{ijk}\theta^i \otimes \theta^j \otimes\theta^k
	\end{align}
    \end{subequations}
with symmetric with respect to the subscripts components $g_{ij}$ and
$h_{ijk}$, we have respectively the following eigenvalues
    \begin{subequations}	\label{11.7}
	\begin{alignat}{3}    \label{11.7a}
&2\lambda_1, &\quad& \lambda_1+\lambda_2 , &\quad& 2\lambda_2
\\    \label{11.7b}
&3\lambda_1 , && 2\lambda_1+\lambda_2, && \lambda_1+2\lambda_2,\quad 3\lambda_3.
	\end{alignat}
    \end{subequations}
	For the quadratic form~\eref{11.6a} hold the implications
	\begin{align}    \label{11.18}
g^\prime = 0
\ &\Rightarrow\
G^{\prime\prime\prime} - 3(\lambda_1 + \lambda_2)G^{\prime\prime}
 + 2(\lambda_1^2 + 4\lambda_1\lambda_2 + \lambda_2^2)G^\prime
 - 4(\lambda_1 + \lambda_2)\lambda_1\lambda_2G = 0
\\    \label{11.19}
G^\prime = 0
\ &\Rightarrow\
g^{\prime\prime\prime} + 3(\lambda_1 + \lambda_2)g^{\prime\prime}
 + 2(\lambda_1^2 + 4\lambda_1\lambda_2
 + \lambda_2^2)g^\prime + 4(\lambda_1 + \lambda_2)\lambda_1\lambda_2g = 0 ,
	\end{align}
which mean that, if the matrix $g:=[g_{ij}]$ is invariant in the flow of $X$,
the form $G$ then satisfies the differential equation in~\eref{11.18}, and
if the form $G$ is invariant in the flow of $X$, its components $g_{ij}$ are
solutions of the r.h.s.\ of~\eref{11.19}.
Now the vector field $\bar{X}$ is
	\begin{equation}    \label{11.23}
\bar{X}
=
- \Bigl(
	\frac{\pd}{\pd g_{11}}, \frac{\pd}{\pd g_{12}}, \frac{\pd}{\pd g_{22}}
   \Bigr)
\cdot
\begin{pmatrix}
      2c_1 & 2c_3 & 0 \\[1ex]
      c_2 & c_1+c_4 & c_3 \\[1ex]
      0 & 2c_2 & 2c_4
\end{pmatrix}
\cdot
\begin{pmatrix} g_{11}\\[1ex] g_{12}\\[1ex] g_{22}\\[1ex] \end{pmatrix} .
	\end{equation}

Similar but more complicated results are valid for the cubic form~\eref{11.6b}.
     \end{Exmp}




    \begin{Exmp}    \label{Exmp11.5}
	For a tensor of type $(1,2)$ (vector-valued quadratic form)
    \begin{equation}    \label{11.9}
K= k_{jl}^{i} R_i\otimes\theta^j\otimes\theta^l
    \end{equation}
with $k_{jl}^{i}=k_{lj}^{i}$, we have the eigenvalues
	\begin{equation}    \label{11.10}
2\lambda_1-\lambda_2 , \quad
\lambda_1,\quad \lambda_2 , \quad
2\lambda_2-\lambda_1 ,
	\end{equation}
where $\lambda_1$ and $\lambda_2$ are double eigenvalues.
The matrix $\bar{C}$ for it is
	\begin{eqnarray}    \label{11.34}
\bar{C}=
- \begin{pmatrix}
      c_1 & 2c_3 & 0 &-c_2 & 0 & 0 \\[1ex]
      c_2 & c_4 & c_3 & 0 & -c_2 & 0 \\[1ex]
      0 & 2c_2 & 2c_4-c_1 & 0 & 0 & -c_2 \\[1ex]
      -c_3 & 0 & 0 & 2c_1-c_4 & 2c_3 & 0 \\[1ex]
      0 & -c_3 & 0 & c_2 & c_1 & c_3 \\[1ex]
      0 & 0 & -c_3 & 0 & 2c_2 & c_4
      \end{pmatrix}
	\end{eqnarray}
and the vector field $\bar{X}$ is (see also~\cite{Boularas&et_al.-2000})
	\begin{align}	\nonumber
\bar{X}=
& -c_1\biggl(k_{11}^1\frac{\pd}{\pd k_{11}^1}
  -k_{22}^1\frac{\pd}{\pd k_{22}^1}+2k_{11}^2\frac{\pd}{\pd k_{11}^2}
  +k_{12}^2\frac{\pd}{\pd k_{12}^2}\biggr)
\\ \nonumber
& -c_2\biggl[-k_{11}^2\frac{\pd}{\pd k_{11}^1}
  +(k_{11}^1-k_{12}^2)\frac{\pd}{\pd h_{12}^1}
  +(2k_{12}^1-k_{22}^2)\frac{\pd}{\pd h_{22}^1}
  +k_{11}^2\frac{\pd}{\pd k_{12}^2}
  +2k_{22}^2\frac{\pd}{\pd k_{22}^2}\biggr]
\\ \nonumber
& -c_3\biggl[2k_{12}^1\frac{\pd}{\pd k_{11}^1}
  +k_{22}^1\frac{\pd}{\pd k_{12}^1}(2k_{12}^2-k_{11}^1)\frac{\pd}{\pd k_{11}^2}
  +(k_{22}^2-k_{12}^1)\frac{\pd}{\pd k_{12}^2}
  -k_{22}^1\frac{\pd}{\pd k_{22}^2}\biggr]
\\    \label{11.37}
& -c_4\biggl(k_{12}^1\frac{\pd}{\pd k_{12}^1}
  +2k_{22}^1\frac{\pd}{\pd k_{22}^1}-k_{11}^2\frac{\pd}{\pd k_{11}^2}
  +k_{22}^2\frac{\pd}{\pd k_{22}^2}\biggr).
	\end{align}
    \end{Exmp}


\section {Conclusion}
\label{Conclusion}

	An invariant of a vector field $X$ on $\field[R]^2$ can be calculated
via the equations
\[
I=\int\omega	\qquad
I=\int \mu\omega \qquad
I=\int\frac{\omega}{\omega(P)}	\qquad
I=u-\int f(v)\Id v
\]
in the cases considered in sections~\ref{Sect0} and~\ref{Sect0-1}. They clearify
the role of infinitesimal symmetry $P$ and the reason why the result can be
expressed via quadratures.

	For a linear vector field, defined via a matrix $C$, we have the
derivative formulae
 $R'=-RC$ and $\theta'=C\theta$
and the exponential $\e^{Ct}$ determines the flow of $X$ and hence the dragging
of tensor fields in that flow. So, the problem for finding the invariants and
flows of linear vector fields is solved completely.

	The considerations presented in this paper admit generalizations in the
multidimensional and infinite dimensional cases.


\section*{Acknowledgments}

	The present paper is done within the Joint research project``Vector
fields and symmetries'' within the bilateral agreement between the Bulgarian
academy of sciences and the Estonian academy of sciences.

	The work of B.~Iliev on this work was partially supported by the
National Science Fund of Bulgaria under Grant No.~F~1515/2005.

	The work of M.\ Rahula was supported by CENS (Tallinn, Estonia).



\addcontentsline{toc}{section}{References}
\bibliography{bozhopub,bozhoref}
\bibliographystyle{unsrt}
\addcontentsline{toc}{subsubsection}{This article ends at page}



\end{document}